\documentclass{article}

\usepackage[latin1]{inputenc}

\usepackage{latexsym}
\usepackage[latin1]{inputenc}
\usepackage{graphicx}
\usepackage{amssymb,amsmath,amsfonts}
\usepackage{enumerate}
\usepackage{amsthm,upgreek}
\usepackage{mathrsfs,yfonts,wasysym,tipa}
\usepackage{stmaryrd}

\newtheorem{teorema}{Theorem}[section]
\newtheorem{prop}[teorema]{Proposition}
\newtheorem{coro}[teorema]{Corollary}
\newtheorem{lema}[teorema]{Lemma}

\begin{document}

\def\R{\mathbb R}
\def\N{\mathbb N}

\title{A discrete approach to Wirtinger's  inequality}
\markright{Wirtinger's  inequality}
\author{Julià Cufí\footnote{The authors were  partially supported by grants 2017SGR358, 2017SGR1725 (Ge\-ne\-ralitat de Catalunya) and MTM2015-66165-P (Ministerio de Economía y Competitividad)}, Agustí Reventós, Carlos J. Rodríguez}

\maketitle


\abstract{Considering  Wirtinger's inequality  for piece-wise equipartite functions we find a discrete version of this classical inequality. The main tool we use is  the theorem of classification of isometries. Our approach provides a new elementary proof of Wirtinger's inequality that also allows to study the case of equality. Moreover it leads  in a natural way to the Fourier series development of $2\pi$-periodic functions.}
\section{Introduction}

The classical Wirtinger inequality states that for a $2\pi$-periodic ${\cal C}^{1}$ function $f(t)$ with $\int_{0}^{2\pi}f(t)\,dt=0$ one has
\begin{eqnarray}\label{W1}\int_{0}^{2\pi}f^{2}(t)\,dt\leq \int_{0}^{2\pi}f'^{2}(t)\,dt,\end{eqnarray}
with equality if and only if $f(t) = a \sin(t) + b \cos(t)$ for some $a,b\in\R$. 

\medskip
The goal of this note is to give a discrete  inequality 
that will imply the above result,  including the case of equality. At the same  time our approach leads in a natural way to the Fourier series development of a $2\pi$-periodic function. 
 
\medskip
Wirtinger did not publish his result, but he communicated it by letter to W. Blascke who included it in \cite{Blaschke16}. The original proof is based on the theory of Fourier series. Discrete
approximations to Wirtinger's inequality have been given by several authors; see for instance \cite{todd}, \cite{shisha}.  

\medskip
As a motivation for a discrete inequality we  consider Wirtinger's inequality for piece-wise equipartite linear functions, that is for continuous functions  $f:[0,2\pi]\longrightarrow \R$,  linear on each interval $[\frac{2\pi}{n}(j-1),\frac{2\pi}{n}j]$, $j=1,\dots, n$ and such that $f(0)=f(2\pi)$. Denoting $f(\frac{2\pi}{n}j)$ by $x_{j}$, $j=1,\dots,n$, and taking $x_{0}=x_{n}$, Wirtinger's inequality for this class of functions is equivalent to 
 the discrete inequality 
 
 \begin{eqnarray}\label{formula150}\sum_{j=1}^{n}x_{j}x_{j-1}\leq \frac{3n^{2}-4\pi^{2}}{3n^{2}+2\pi^{2}}\, ,\end{eqnarray}
for $x_{j}\in\R$, $j=1,\dots,n$, $x_{0}=x_{n}$, $\sum_{j=1}^{n}x_{j}=0$  and  $\sum_{j=1}^{n}x_{j}^{2}=1$. 

\medskip

Wirtinger's inequality can then be obtained from the above inequality by a limiting process.  

\medskip
We shall obtain \eqref{formula150} as a consequence of the following 

\medskip
{\bf Theorem 2.1.} {\em Let $x_{1},\dots,x_{n}\in \R$, for $n\geq 4$,  with $\sum_{i=1}^{n}x_{i}=0$ and $\sum_{i=1}^{n}x_{i}^{2}=1$. Then
\begin{eqnarray}\label{W3}\sum_{i=1}^{n}x_{i}x_{i-1}\leq \cos(\frac{2\pi}{n}), \end{eqnarray}
with $x_{0}=x_{n}$.
Equality holds if and only if $$x_{i}=a \cos(\frac{2\pi}{n}i)+ b \sin(\frac{2\pi}{n}i), \quad i=1,\dots,n, $$
for   $ a, b\in\R  $ satisfying  $a^{2}+b^{2}=2/n$.}

\medskip
This result that can be considered as the {\em Wirtinger discrete  inequality} was obtained by Fan, Taussky and Todd in \cite{todd} where it is used  to obtain classical Wirtinger inequality \eqref{W1} but,  as the authors say, without the equality clause. Other proofs of Theorem \ref{teorema21} have been published later, see for instance \cite{shisha}. 

For completeness 
we provide here a simple different  proof 
of the above  result based on  
the theorem of classification of isometries applied to the cyclic isometry $T$ given by $$T(x_{1},x_{2},\dots,x_{n})=(x_{n}, x_{1}, x_{2}, \dots, x_{n-1}), $$
since the left hand-side of \eqref{W3} can be written as  $\langle X, T(X)\rangle$, where $X=(x_{1},x_{2},\dots,x_{n})$.

As we have said our approach,  based on inequality \eqref{formula150}, leads to inequality \eqref{W1}, and allows to caracterize functions for which equality holds.  This characterization is somewhat delicate but the argument used has 
a surprising consequence: the   Fourier series development of a $2\pi$-periodic function.


\section{Discrete Wirtinger's inequality}\label{preli}
In order to find a discrete version of the Wirtinger inequality
we consider this inequality for piece-wise equipartite linear functions.

For $n\in\N$, $n\geq 4$,  let $f:[0,2\pi]\longrightarrow \R$ be a continuous function, linear on each interval $[\frac{2\pi}{n}(j-1),\frac{2\pi}{n}j]$, $j=1,\dots, n$ and such that $f(0)=f(2\pi)$. 
Denoting $f(\frac{2\pi}{n}j)$ by $x_{j}$, $j=1,\dots,n$, and taking $x_{0}=x_{n}$, a computation shows that 

\begin{eqnarray}\label{formula0}\int_{0}^{2\pi}f^{2}(t)\,dt=\frac{2\pi}{3n}\sum_{j=1}^{n}(2x_{j}^{2}+x_{j}x_{j-1})\end{eqnarray}
and
\begin{eqnarray}\label{formula1}\int_{0}^{2\pi}f'^{2}(t)\,dt=\frac{n}{\pi}\sum_{j=1}^{n}(x_{j}^{2}-x_{j}x_{j-1}).\end{eqnarray}
So the inequality $$\int_{0}^{2\pi}f^{2}(t)\,dt\leq \int_{0}^{2\pi}f'^{2}(t)\,dt$$
is equivalent to
\begin{eqnarray}\label{formula9}\sum_{j=1}^{n}x_{j}x_{j-1}\leq \frac{3n^{2}-4\pi^{2}}{3n^{2}+2\pi^{2}}\sum_{j=1}^{n}x_{j}^{2}.\end{eqnarray}

Assuming now $\int_{0}^{2\pi}f(t)\,dt=0$, that means $\sum_{i=1}^{n}x_{j}=0$, it follows that  Wirtinger's inequality for piece-wise linear functions
is equivalent to  \eqref{formula9} with this  additional hypothesis or, normalizing, 

\begin{eqnarray}\label{formula15}\sum_{j=1}^{n}x_{j}x_{j-1}\leq \frac{3n^{2}-4\pi^{2}}{3n^{2}+2\pi^{2}},\quad \mbox{with }\sum_{j=1}^{n}x_{j}=0, \; \sum_{j=1}^{n}x_{j}^{2}=1.\end{eqnarray}

This is a problem of maximizing a given quadratic form under some restrictions.
It can be solved by different methods such as Lagrange multipliers or by the determination of the least characteristic value of a Hermitian matrix, as done in \cite{todd}. As said our approach is based on the theorem of classification of isometries. 

\subsection*{The canonical expression of the quadratic form.}

The left-hand side of \eqref{formula15} leads in a natural way to consider the  cyclic isometry
$$T(x_{1},x_{2},\dots,x_{n})=(x_{n}, x_{1}, x_{2}, \dots, x_{n-1}),$$
since
$$\sum_{j=1}^{n}x_{j}x_{j-1}=\langle X,T(X)\rangle,$$ where $X=(x_{1},\dots,x_{n})$
and $\langle \,,\, \rangle$ is the standard scalar product.
Hence, in order to prove \eqref{formula15} we start by analyzing the structure of the isometry $T$. This will allow us to find the canonical expression of the quadratic form $\langle X, T(X)\rangle$.


The theorem of classification of isometries  (see \cite{agusti}) applied to $T$ asserts that there is an orthonormal basis $(e_{1},\dots,e_{n})$ such that, denoting $\alpha_{k}=\frac{2\pi}{n}k$,  one has  for $n$ even 
\begin{eqnarray*}
T(e_{1})&=&e_{1},\\
T(e_{2})&=&-e_{2},\\
T(e_{2k+1})&=&(\cos \alpha_{k})\,e_{2k+1}+(\sin \alpha_{k})\,e_{2k+2},\\
T(e_{2k+2})&=&(-\sin \alpha_{k})\,e_{2k+1}+(\cos \alpha_{k})\,e_{2k+2}, \qquad k=1,\dots, (n-2)/2,
\end{eqnarray*}
and for $n$ odd 

\begin{eqnarray*}
T(e_{1})&=&e_{1},\\
T(e_{2k})&=&(\cos \alpha_{k})\,e_{2k}+(\sin \alpha_{k})\,e_{2k+1},\\
T(e_{2k+1})&=&(-\sin \alpha_{k})\,e_{2k}+(\cos \alpha_{k})\,e_{2k+1}, \qquad k=1,\dots, (n-1)/2.
\end{eqnarray*}
In fact, it can be seen by using elementary trigonometric formulas that for $n$ even, this basis is given by
\begin{eqnarray}\label{odd20}
e_{1}&=&\frac{1}{\sqrt{n}}(1,\dots,1),\nonumber\\
e_{2}&=&\frac{1}{\sqrt{n}}(1,-1\dots,1,-1),\nonumber\\
e_{2k+1}&=&\sqrt{\frac{2}{n}}\left(1,\cos(\frac{2\pi }{n}k), \cos(\frac{2\pi }{n}2k),\dots, \cos(\frac{2\pi }{n}(n-1)k)\right),\nonumber\\\\
e_{2k+2}&=&\sqrt{\frac{2}{n}}
\left(0,\sin(\frac{2\pi }{n}k), \sin(\frac{2\pi }{n}2k),\dots, \sin(\frac{2\pi }{n}(n-1)k)\right),\nonumber\end{eqnarray}
and for $n$ odd by
\begin{eqnarray}\label{odd2}
e_{1}&=&\frac{1}{\sqrt{n}}(1,\dots,1),\nonumber\\
e_{2k}&=&\sqrt{\frac{2}{n}}\left(1,\cos(\frac{2\pi }{n}k), \cos(\frac{2\pi }{n}2k),\dots, \cos(\frac{2\pi }{n}(n-1)k)\right),\nonumber\\\\
e_{2k+1}&=&\sqrt{\frac{2}{n}}
\left(0,\sin(\frac{2\pi }{n}k), \sin(\frac{2\pi }{n}2k),\dots, \sin(\frac{2\pi }{n}(n-1)k)\right)\cdot \nonumber\end{eqnarray}
Since  $\langle e_{i}, T(e_{j})\rangle+\langle e_{j}, T(e_{i})\rangle=0, \mbox{ for } i\neq j, $  we get 
 for every vector $X=\sum_{i=1}^{n}y_{i}e_{i}$,
 
 \begin{eqnarray*}
 \langle X,T(X)\rangle&=&\sum_{i,j=1}^{n}y_{i}y_{j}\langle e_{i}, T(e_{j})\rangle= \sum_{i=1}^{n}y_{i}^{2}\langle e_{i}, T(e_{i})\rangle.
 \end{eqnarray*}
 Hence the canonical expression of the quadratic form $\langle X, T(X)\rangle$ is  for even  $n$  \begin{eqnarray}\label{even}\langle X,T(X)\rangle=y_{1}^{2}-y_{2}^{2}+\sum_{k=1}^{(n-2)/2}(y_{2k+1}^{2}+y_{2k+2}^{2})\cos \alpha_{k},\end{eqnarray}  and  for odd $n$ 
 \begin{eqnarray}\label{odd}\langle X,T(X)\rangle=y_{1}^{2}+\sum_{k=1}^{(n-1)/2}(y_{2k}^{2}+y_{2k+1}^{2})\cos \alpha_{k}.\end{eqnarray}

\subsection*{The discrete inequality}
The maximum of the quadratic form $\langle X, T(X)\rangle$ is given by the following

\begin{teorema}[Discrete Wirtinger's inequality]\label{teorema21} Let $x_{1},\dots,x_{n}\in \R$, for $n\geq 4$, with $\sum_{i=1}^{n}x_{i}=0$ and $\sum_{i=1}^{n}x_{i}^{2}=1$. Then
$$\sum_{i=1}^{n}x_{i}x_{i-1}\leq \cos(\frac{2\pi}{n}), $$
with $x_{0}=x_{n}$.
Equality holds if and only if \begin{eqnarray*}x_{i}=a \cos(\frac{2\pi}{n}i)+ b \sin(\frac{2\pi}{n}i), \quad i=1,\dots,n, \end{eqnarray*}
for  $a, b\in\R$ satisfying  $a^{2}+b^{2}=2/n$.
\end{teorema}
{\em Proof.}  With the previous notation we must prove 
$$\langle X, T(X)\rangle \leq \cos(\frac{2\pi}{n}).$$

%
%

Since $ \langle X, e_{1}\rangle=0$, $||X||=1$ it is  $X=\sum_{j=2}^{n}y_{j}e_{j}$,  $\sum_{j=2}^{n}y_{j}^{2}=1$  and we get from \eqref{even}  
 
 $$\langle X,T(X)\rangle \leq \cos(\frac{2\pi}{n})\sum_{k=1}^{(n-2)/2}(y_{2k+1}^{2}+y_{2k+2}^{2})\leq \cos(\frac{2\pi}{n})$$
  for $n$ even, and from \eqref{odd}
 $$\langle X,T(X)\rangle \leq \cos(\frac{2\pi}{n})\sum_{k=1}^{(n-1)/2}(y_{2k}^{2}+y_{2k+1}^{2})=\cos(\frac{2\pi}{n})$$
 for $n$ odd. This proves the first part of the Lemma.
 
 Equality holds when  $X=y_{3}e_{3}+y_{4}e_{4}$ for $n$ even and $X=y_{2}e_{2}+y_{3}e_{3}$ for $n$ odd. Substituting $e_{2},e_{3}, e_{4}$ by the expressions in \eqref{odd20} and \eqref{odd2} the Lemma follows. $\square$

 \medskip 
 As a consequence of this result we obtain inequality \eqref{formula15}.

%
 
 \begin{prop}\label{WD} Let $x_{1},\dots,x_{n}\in \R$, for $n\geq 4$,  with $\sum_{i=1}^{n}x_{i}=0$ and $\sum_{i=1}^{n}x_{i}^{2}=1$. Then \begin{eqnarray}\label{etiqueta}\sum_{i=1}^{n}x_{i}x_{i-1}<  \frac{3n^{2}-4\pi^{2}}{3n^{2}+2\pi^{2}},\end{eqnarray} with $ x_{0}=x_{n}$. 
 \end{prop}
{\em Proof.} By Theorem  \ref{teorema21} in order to prove \eqref{etiqueta} it is enough to show that $$\cos(\frac{2\pi}{n})< \frac{3n^{2}-4\pi^{2}}{3n^{2}+2\pi^{2}}.$$

Denoting $2\pi/n$ by $\alpha$ the above inequality is equivalent to
$$\cos\alpha< \frac{6-2\alpha^{2}}{6+\alpha^{2}},$$
which  using that
$\cos\alpha < 1-\alpha^{2}/2+\alpha^{4}/24$ is easily verified. 
 \qed
 
  We remark that equality in \eqref{formula15}  never holds.

 \begin{coro}\label{prop222} For $n\in \N$, $n\geq 4$, let $f:[0,2\pi]\longrightarrow \R$ be a continuous function, linear on each interval $[\frac{2\pi}{n}(j-1),\frac{2\pi}{n}j]$, $j=1,\dots, n$ and such that $f(0)=f(2\pi)$.  Assume that 
 $\int_{0}^{2\pi}f(t)dt=0$.
 Then
 \begin{eqnarray}\label{prop22}\int_{0}^{2\pi} f^{2}(t)\,dt\leq \int_{0}^{2\pi}f'^{2}(t)\,dt.\end{eqnarray}
\end{coro}

{\em Proof.} As said, inequality \eqref{prop22} with hypothesis $\int_{0}^{2\pi}f(t)\,dt=0$ is equivalent to \eqref{formula15}. So the Corollary is a direct consequence of Proposition  \ref{WD}. 
 $\square$

\medskip

\section{Wirtinger's inequality}
Now we can obtain, by a limiting process,  the classical Wirtinger's inequality.

\begin{teorema}[Wirtinger's inequality]\label{T23}Let $f:\R\longrightarrow\R$ be a $2\pi$-periodic ${\cal C}^{1}$ function such that $\int_{0}^{2\pi}f(t)\,dt=0.$
Then 
\begin{eqnarray}\label{ineq}\int_{0}^{2\pi}f^{2}(t)dt\leq \int_{0}^{2\pi}f'^{2}(t)\,dt.\end{eqnarray}
Equality holds if and only if $f(t)=a\cos(t)+b\sin(t)$ for some $a,b\in\R$.
\end{teorema}
{\em Proof.} For each  $n\in \N$, $n\geq 4$, let $\phi_{n}(t)$ be the function  linear on each interval $[\frac{2\pi}{n}(j-1),\frac{2\pi}{n}j]$,  with $\phi_{n}(\frac{2\pi}{n}j)=f(\frac{2\pi}{n}j)$, $j=1,\dots, n$.

Set $x_{j,n}=f(\frac{2\pi}{n}j)$, $m_{n}=\dfrac{1}{n}\sum_{j=1}^{n}x_{j}$, and $\tilde{x}_{j,n}=x_{j,n}-m_{n}$. Let $\tilde{\phi}_{n}(t)$ be the function  linear on each interval $[\frac{2\pi}{n}(j-1),\frac{2\pi}{n}j]$,  with $\tilde{\phi}_{n}(\frac{2\pi}{n}j)=\tilde{x}_{j,n}$, $j=1,\dots, n$. Equivalently, $\tilde{\phi}_{n}(t)=\phi_{n}(t)-m_{n}$.

Since $\int_{0}^{2\pi}\tilde{\phi}(t)\,dt=0$ it follows, by Corollary \ref{prop222},  that $$\int_{0}^{2\pi}\tilde{\phi}_{n}^{2}(t)\,dt\leq \int_{0}^{2\pi}\tilde{\phi}_{n}'^{2}(t)\, dt.$$
Moreover since $f$ is a ${\cal C}^{1}$ function we have
$$\lim_{n\to \infty}\int_{0}^{2\pi}\phi_{n}^{2}(t)\,dt=\int_{0}^{2\pi}f^{2}(t)\,dt$$
and 
\begin{eqnarray}\label{formula2}\lim_{n\to \infty}\int_{0}^{2\pi}\phi_{n}'^{2}(t)\,dt=\int_{0}^{2\pi}f'^{2}(t)\,dt.\end{eqnarray}
Finally the hypothesis $\int_{0}^{2\pi}f(t)\,dt=0$ yields $\lim_{n\to\infty}m_{n}=0$ and so
\begin{eqnarray*}\lim_{n\to\infty}\int_{0}^{2\pi}\tilde{\phi}_{n}^{2}(t)\,dt=	\lim_{n\to\infty}\int_{0}^{2\pi}\phi_{n}^{2}(t)\,dt,\end{eqnarray*}
and inequality \eqref{ineq} follows. 

\medskip
It remains to analize when equality holds in \eqref{ineq}. 

From now on we will assume that $n$ is an  odd integer; the case $n$ even is dealt similarly. Let $H_{k}=\langle e_{2k}, e_{2k+1}\rangle$ denote the subspace of $\R^{n}$ generated by $e_{2k}$ and $e_{2k+1}$, the vectors introduced in Section \ref{preli}, for $k=1,\dots, (n-1)/2$.  Let $P_{k}$ be the orthogonal projection from $\R^{n}$ on $H_{k}$, and let $P_{0}$ be the orthogonal projection on $H_{0}=\langle e_{1}\rangle$.

\begin{lema}\label{lema}Let  $f$ be a function satisfying  the hypotheses of Theorem \ref{T23} and such that  equality holds in \eqref{ineq}.
For each $n\geq 4$ let $X_{n}$ be the vector of components $x_{j,n}=f(\frac{2\pi}{n}j)$, $j=1,\dots, n$. Then

$$\lim_{n\to\infty}\frac{1}{n}\sum_{k=2}^{(n-1)/2}||P_{k}(X_{n})||^{2}= 0.$$  

\end{lema}
{\em Proof}.  By the definition of Riemann's integral we have
$$\lim_{n\to\infty}\frac{2\pi}{n}||X_{n}||^{2}=\lim_{n\to\infty}\frac{2\pi}{n}\sum_{k=0}^{p}||P_{k}(X_{n})||^{2}=\int_{0}^{2\pi}f^{2}(t)\,dt,$$ where $p=(n-1)/2$.

From \eqref{formula1}, \eqref{odd} and \eqref{formula2}
it follows that 
\begin{eqnarray*}&&\lim_{n\to\infty}\frac{n}{\pi}\left(\sum_{k=0}^{p}||P_{k}(X_{n})||^{2}-\langle T(X_{n}),X_{n}\rangle \right)\\&=&\lim_{n\to\infty}\frac{n}{\pi}\sum_{k=1}^{p}(||P_{k}(X_{n})||^{2}(1-\cos(\frac{2\pi}{n}k)))=\int_{0}^{2\pi}f'^{2}(t)\,dt.\end{eqnarray*}

As a consequence of equality in \eqref{ineq} we get
$$\lim_{n\to \infty}\sum_{k=1}^{p}\left[\frac{n}{\pi}(1-\cos(\frac{2\pi}{n}k))-\frac{2\pi}{n}\right]||P_{k}(X_{n})||^{2}=0.$$

The Lemma follows from the inequality

$$\frac{n}{\pi}(1-\cos(\frac{2\pi}{n}k))-\frac{2\pi}{n}\geq \frac{1}{n}$$ which is true for $k\geq 2$ (which implies $n\geq 5$) using that $\cos(x)\leq 1-x^{2}/2+x^{4}/24$. $\square$

\medskip 
To continue the proof of Theorem \ref{T23}, for each vector $X=(x_{1},\dots,x_{n})$ let $L_{X}$ be the function that is linear on each interval $[\frac{2\pi}{n}(j-1),\frac{2\pi}{n}j]$ with $L_{X}(\frac{2\pi}{n}j)=x_{j}$, $j=1,\dots,n$, ($x_{0}=x_{n}$).

When $X_{n}$ is the vector of components $x_{j,n}=f(\frac{2\pi}{n}j)$, $j=1,\dots, n$, $L_{X_{n}}$ is the function $\phi_{n}$ defined at the begining  of this  proof. So we can assume that $\sum_{j=1}^{n}x_{j,n}=0$ and we know that $\lim_{n\to \infty}L_{X_{n}}=f$.

Writing $X_{n}=y_{2}e_{2}+y_{3}e_{3}+\sum_{k=2}^{p}(y_{2k}e_{2k}+y_{2k+1}e_{2k+1})$
we have
$$L_{X_{n}}=y_{2}L_{e_{2}}+y_{3}L_{e_{3}}+\sum_{k=2}^{p}(y_{2k}L_{e_{2k}}+y_{2k+1}L_{e_{2k+1}}):=\alpha_{n}+\beta_{n}.$$

To finish the proof we need to show that
$$\lim_{n\to\infty}\alpha_{n}=a\cos(t)+b\sin(t), \mbox{ for some $a,b\in\R$ and } \lim_{n\to\infty}\beta_{n}=0.$$ 

Formula \eqref{formula0} can be writen as
$$\int_{0}^{2\pi}L_{X}^{2}dt=\frac{4\pi}{3n}||X||^{2}+\frac{2\pi}{3n}\langle X, T(X)\rangle$$
which gives, by using the identity of polarization, 
$$\langle L_{X},L_{Y}\rangle:=\int_{0}^{2\pi}L_{X}L_{Y}dt=\frac{4\pi}{3n}\langle X,Y\rangle+\frac{\pi}{3n}\langle T(X), Y\rangle+\frac{\pi}{3n}\langle X, T(Y)\rangle$$ for  two vectors $X,Y$.

In particular one gets $\langle L_{e_{i}},L_{e_{j}}\rangle=0$, $i\neq j$, and hence 
$$\langle L_{X},L_{e_{j}}\rangle=y_{j}\langle L_{e_{j}},L_{e_{j}}\rangle,\quad j=2,\dots,n,$$ and

$$\langle L_{e_{2k}},L_{e_{2k}}\rangle=\frac{4\pi}{3n}(1+\frac{1}{2}\cos\frac{2\pi}{n}k)=\langle L_{e_{2k+1}},L_{e_{2k+1}}\rangle, \quad k=1,\dots,\frac{n-1}{2}.$$

Writing $\tilde{e}_{2}=\sqrt{\frac{n}{2}}\,e_{2}$, $\tilde{e}_{3}=\sqrt{\frac{n}{2}}\,e_{3}$
we have

\begin{eqnarray*}
\lim_{n\to\infty}\alpha_{n}&=&\lim_{n\to\infty}(y_{2}L_{e_{2}}+y_{3}L_{e_{3}})\\&=&\lim_{n\to\infty}\frac{\langle L_{X_{n}}, L_{\tilde{e}_{2}}\rangle \frac{\sqrt{2}}{\sqrt{n}}L_{\tilde{e}_{2}}\frac{\sqrt{2}}{\sqrt{n}}+\langle L_{X_{n}}, L_{\tilde{e}_{3}}\rangle \frac{\sqrt{2}}{\sqrt{n}}L_{\tilde{e}_{3}}\frac{\sqrt{2}}{\sqrt{n}}}{\frac{4\pi}{3n}(1+\frac{1}{2}\cos\frac{2\pi}{n})}
\\&=&\frac{1}{\pi}\left(\int_{0}^{2\pi}f(t)\cos(t)\,dt\right) \cos t+\frac{1}{\pi}\left(\int_{0}^{2\pi}f(t)\sin(t)\,dt\right) \sin t.
\end{eqnarray*}

Thus $\lim_{n\to\infty}\alpha_{n}=a\cos t+ b\sin t$, as wanted, where $a,b$ are the first Fourier coefficients of $f$. 

As concerning $\lim_{n\to\infty}\beta_{n}$ we have
\begin{eqnarray*}\langle \beta_{n}, \beta_{n}\rangle&=&\int_{0}^{2\pi}\beta_{n}\cdot \beta_{n}\, dt=\\&=&\sum_{k=2}^{p}(y_{2k}^{2}+y_{2k+1}^{2})\frac{4\pi}{3n}(1+\frac{1}{2}\cos\frac{2\pi}{n}k)\leq 2\pi\frac{1}{n}\sum_{k=2}^{p}||P_{k}(X_{n})||^{2},\end{eqnarray*}
and the proof finishes by applying Lema \ref{lema}. $\square$

\bigskip

{\bf Remark.}
%
%
Let $f$ be a $2\pi$-periodic ${\cal C}^{1}$ function such that $\int_{0}^{2\pi}f(t)\,dt=0$. The same argument used to calculate $\lim_{n\to\infty}\alpha_{n}$ in the above proof, applied also to $\beta_{n}$
shows that $f$ can be written as $$ f(t)=\sum_{j=1}^{\infty}(a_{j}\cos(jt)+b_{j}\sin(jt)),$$ with $$a_{j}=\frac{1}{\pi}\int_{0}^{2\pi}f(t)\cos(jt)dt,\quad b_{j}=\frac{1}{\pi}\int_{0}^{2\pi}f(t)\sin(jt)dt.$$ 
If we drop the assumption $\int_{0}^{2\pi}f(t)\,dt=0$ we need to add in the above  expression of $f$ the term $\frac{1}{2\pi}\int_{0}^{2\pi}f(t)\,dt$. So, the discrete approach we have developped here leads, in a natural way, to the  well known Fourier series development  of a $2\pi$-periodic function. 
\vspace{1cm}
\bibliographystyle{plain}
\bibliography{Bibliografia}

\noindent {\em Departament de Matemàtiques, Universitat Aut\`{o}noma de Barcelona\\ 08193 Bellaterra, Barcelona, Catalonia\\

\noindent jcufi@mat.uab.cat,  agusti@mat.uab.cat, rdrgzcarlos@gmail.com.

\end{document}